\documentclass[a4paper]{article}
\usepackage[top=25truemm,bottom=25truemm,left=20truemm,right=20truemm]{geometry}
\usepackage{amsmath,amsthm,amsfonts,amssymb,stmaryrd,indentfirst,tabularray,makecell,hhline,comment,mathtools}
\usepackage[mathscr]{euscript}
\usepackage[shortlabels]{enumitem}
\usepackage[labelfont=rm]{caption}
\usepackage{titlesec}
\usepackage{tikz-cd}
\usepackage{url}
\usepackage[backref]{hyperref}
\hypersetup{
    colorlinks = true,
    linkcolor = [rgb]{0.3,0.55,0.25},
    citecolor = [rgb]{0.164,0.39,0.77},
    urlcolor = [rgb]{0.7,0.2,0.4}
}

\usepackage{tikz}
\usetikzlibrary{decorations.markings}

\newcommand{\RGbar}{
	\begin{tikzpicture}[line width=0.6pt, baseline=-0.65ex, scale=1.2]
	
	\node[fill, circle, inner sep=0pt, minimum size=3pt] at (0,0) (n1) {};
	\node[fill, circle, inner sep=0pt, minimum size=3pt] at (0.5,0) (n2) {};
	
	\draw[-Stealth] (n1)--(n2);
	
	\end{tikzpicture}
}

\newcommand{\RGtadpole}{
	\begin{tikzpicture}[line width=0.6pt, baseline=-1ex, scale=1.2]
	
	\node[fill, circle, inner sep=0pt, minimum size=3pt] at (0,0) (n1) {};
	
	\draw (n1) to[out=90, in=0] (-0.2,0.22) to[out=180, in=90] (-0.4,0) to[out=-90, in=180] (-0.2,-0.2);
	\draw[-Stealth] (-0.2,-0.2) to[out=0, in=238] (n1);
	
	\end{tikzpicture}
}

\DeclareMathOperator{\id}{id}
\DeclareMathOperator{\Ker}{Ker}
\DeclareMathOperator{\Imag}{Im}
\DeclareMathOperator{\Aut}{Aut}
\DeclareMathOperator{\Der}{Der}
\DeclareMathOperator{\Tr}{Tr}

\newcommand{\gr}[0]{_\mathrm{gr}}

\theoremstyle{plain}
\newtheorem{theorem}{Theorem}
\newtheorem*{theorem*}{Theorem}

\theoremstyle{definition}

\newtheorem{remark}[theorem]{Remark}
\newtheorem{example}[theorem]{Example}

\newtheorem{question}[theorem]{Question}

\numberwithin{theorem}{section}
\everymath{\displaystyle}
\allowdisplaybreaks

\usepackage{apptools}
\titleformat{\section}{\large\scshape}{\IfAppendix{\appendixname}{} \thesection.}{3pt}{}
\titleformat{\subsection}{\scshape}{\thesubsection.}{3pt}{}

\usepackage{titling}  
\pretitle{
	\begin{center}\Large\bfseries
}
\posttitle{
	\end{center}\ \\[-20pt]
}
\preauthor{
	\begin{center}\large
}
\postauthor{
	\end{center}\ \\[-60pt]
}

\title{One-dimensionality of Certain Cocycles\\for the Detection of the Johnson Cokernel}
\author{Toyo TANIGUCHI \thanks{Graduate School of Mathematical Sciences, The University of Tokyo. 3-8-1, Komaba, Meguro-ku, Tokyo, 153-8914, Japan. E-mail: \texttt{toyo(at)ms.u-tokyo.ac.jp}}}
\date{}

\begin{document}
\maketitle

\centerline{\it Dedicated to Professor Shigeyuki Morita on the occasion of his 80th birthday}

\begin{abstract}
The Enomoto--Satoh (ES) trace detects the Johnson cokernel, and its $1$-cocycle property is important for the proof that the Johnson image is annihilated by the ES trace. Via the natural map from the ribbon graph complex introduced by Merkulov and Willwacher to the Chevalley--Eilenberg complex of the Lie algebra of symplectic derivations, where the Johnson image lives, the ES trace is essentially obtained from the $1$-cocycle given by the unique ribbon graph with one vertex and one edge. In this perspective, this ribbon graph is the ``universal'' version of the ES trace. The main result of this paper is that there are no other (linearly independent) $1$-cocycles in the ribbon graph complex, showing that nothing can be found there for the detection of the Johnson cokernel. The proof is done by applying the result of Church--Farb--Putman and Morita--Sakasai--Suzuki, which states that the virtual top-dimensional cohomology of the moduli space of marked Riemann surfaces vanishes.\\
\end{abstract}

\noindent{\textit{2020 Mathematics Subject Classification: 18M70, 20F34, 32G15.}

\section{Introduction}
Let $\mathbb{K}$ be a field of characteristic zero, $\Sigma = \Sigma_{g,1}$ a connected compact oriented surface of genus $g$ with one boundary component, and $H = H_1(\Sigma;\mathbb{K})$ its first homology group. In the study of the mapping class group $\mathcal M$ of the surface $\Sigma$, Johnson introduced a descending filtration $\{\mathcal M(k)\}_{k\geq 1}$ on $\mathcal M$. The \textit{Johnson homomorphism} (in the formulation of Morita's refined version \cite{morita} based on the original one in \cite{johnson1} and \cite{johnson2}) is an injective map of Lie algebras
\[
	\tau\colon \mathrm{gr}(\mathcal{M}(1)\otimes_\mathbb{Z}\mathbb{K}) \to \Der_\omega^+(L(H))
\]
to the pro-nilpotent Lie algebra of symplectic derivations on the free Lie algebra $L(H)$ over $H$. This allows us to embed (the associated graded of) the \textit{Torelli group} $\mathcal M(1)$ into an algebraically more tractable space. Furthermore, for $g \geq 3$, Hain showed in \cite{hain1} and \cite{hain2} that its image is generated by the degree-1 subspace of $\Der_\omega^+(L(H))$ and gave an explicit presentation for $g\geq 6$.

Morita's problem concerning the Johnson homomorphism was to determine the cokernel of $\tau$ as an $\mathrm{Sp}(H)$-representation (see Problem 32 of \cite{morita2}), and he introduced the \textit{Morita trace} $\Tr^\mathsf{M}$ in \cite{morita} for the detection (or as an ``invariant'') of the cokernel: the map $\Tr^\mathsf{M}$ is non-trivial and $\Tr^\mathsf{M}\circ \,\tau = 0$ holds on the degree $\geq 2$ part. Later, more stronger \textit{Enomoto--Satoh trace} $\Tr^\mathsf{ES}$ was constructed in \cite{es}, which also satisfies $\Tr^\mathsf{ES}\circ\,\tau = 0$ on the degree $\geq 2$ part, and further enhanced by Conant, Kassabov and Vogtmann based on their series of work \cite{ckv1}, \cite{ckv2}, \cite{conant} and \cite{ck}. As for recent progress related to the Torelli group and the Johnson homomorphism (in the case of $\Sigma_{g,1}$ or the closed surface $\Sigma_g$), Morita--Sakasai--Suzuki \cite{mss2} determined the image and cokernel up to degree $6$ based on computer calculation, Garoufalidis--Getzler \cite{garget} gave an explicit formula for the $\mathfrak{sp}$-character of the pro-unipotent completion of the Torelli group of $\Sigma_g$ in the stable range under the Koszulity assumption, and Kupers--Randal-Williams \cite{krw} showed that the kernel of the Johnson homomorphism stably comprises of only the trivial representation. The Koszulity was also proved by both Felder--Naef--Willwacher \cite{fnw} and Kupers--Randal-Williams \cite{krw}.

 The ES trace is shown to be equivalent to the graded version of the Turaev cobracket
 \[
 	\delta\gr\colon |T(H)| \to |T(H)|\otimes |T(H)|
\]
where $|T(H)| = T(H)/[T(H),T(H)]$ is the \textit{trace space} of the tensor algebra over $H$ (see \cite{akkn}, Section 9). For more details, see the survey \cite{gtjohnson} by Kawazumi and Kuno. On the other hand, the map $\delta\gr$ naturally appears from the ribbon graph complex in the work \cite{ribbon} of Merkulov and Willwacher. They constructed the properad $\mathcal{RG}ra_d$ of (oriented and labelled) ribbon graphs for $d\in\mathbb{Z}$, and showed that
there is the canonical representation
\[
	\rho_H \colon \mathcal{RG}ra_1 \to \mathcal End_{|T(H)|}
\]
under which $G_1 = \RGtadpole$, the unique ribbon graph with one vertex and one edge, is sent to the graded version of the Turaev cobracket $\delta\gr$. We also have the morphism of properads
\[
	s_1^*\colon \mathcal Lieb_d \to \mathcal{RG}ra_d
\]
from the properad governing degree-shifted Lie bialgebras, which defines the \textit{ribbon graph complex} $\mathsf{RGC}_d$ using the framework of the deformation complex in \cite{merval}. Its cohomology contains the direct product of the cohomology of the moduli spaces of surfaces for all genera by the result of Kontsevich; see \eqref{eq:moduliisom} for a detailed statement.

The key property of the ES trace is that it is a Lie algebra $1$-cocycle on $\Der^+_\omega(L(H))$, which is regarded as a subspace of $|T(H)|$ in a canonical way. Since the ES trace is obtained from a $1$-cocycle $G_1$ in $\mathsf{RGC}_1$ in the form of $\delta\gr$, we want to find more $1$-cocycles for the detection of the Johnson cokernel. However, the main result of this paper is the following:
\begin{theorem*}[Theorem \ref{thm:nococyc}]
The space of $1$-cocycles  is $1$-dimensional and it is spanned by $G_1$.
\end{theorem*}
The proof is based on the remarkable result of Church--Farb--Putman \cite{cfp} and Morita--Sakasai--Suzuki \cite{mss} that states the virtual top-dimensional cohomology of $H^{4g-3}(\mathbf{M}_g^1;\mathbb{Q})$ of the moduli space of Riemann surfaces of genus $g$ and one marked point vanishes.\\

The construction of the trace map by Conant, Kassabov and Vogtmann is based on Lie graphs, and it takes values in a module $\Omega_n$ which surjects onto the cohomology $H^{2n-3}(\mathrm{Out}(F_n); \overline{T(V)^{\otimes n}})$ of the outer automorphism group of the free group with a certain coefficient module. The discussion in this paper is based on ribbon graphs ($=$ ``associative'' graphs), but since the Johnson image naturally lives in the space of derivations of a free Lie algebra, it would be interesting if we can do similar approach in the Lie world (see Question \ref{q:Lie}) for further understanding of the Johnson cokernel.\\

\noindent\textbf{Organisation of the paper.} Section \ref{sec:Johnson} is a review of the Johnson homomorphism, and Section \ref{sec:rg} recalls the ribbon graph complex in the formulation by Merkulov and Willwacher. The main result and its proof occupy Section \ref{sec:mainthm}.\\

\noindent\textbf{Acknowledgements.} The author thanks Sergei Merkulov for kindly explaining the isomorphism between the cohomology of the ribbon graph complex and the moduli space of Riemann surfaces, Takuya Sakasai for pointing out some references on this topic, and Nariya Kawazumi for careful reading of the draft. This work was supported by JSPS KAKENHI Grant Number 25KJ0734.\\

\noindent\textbf{Conventions.} $\mathbb{K}$ is a field of characteristic zero. Unadorned tensor products are always over $\mathbb{K}$.\\

\section{The Johnson Homomorphism and the Enomoto--Satoh Trace}\label{sec:Johnson}

In this section, we recall the Johnson homomorphism and the Enomoto--Satoh trace. The material here is contained in the survey \cite{gtjohnson} by Kawazumi and Kuno and Section 9 of \cite{akkn}. For the original works by Johnson, see \cite{johnson1} and \cite{johnson2}.\\

Let $\Sigma = \Sigma_{g,1}$ be a connected oriented compact surface of genus $g>0$ with one boundary component, $\pi = \pi_1(\Sigma, *)$ the fundamntal group of $\Sigma$ with the base point $*\in\partial\Sigma$, and $H_\mathbb{Z} = H_1(\Sigma;\mathbb{Z})$ the first homology group. The lower central series $\{\Gamma_k\} _{k\geq 1}$ of $\pi$ recursively given by 
\[
	\Gamma_1 = \pi\;\mbox{ and }\; \Gamma_{k+1} = [\pi, \Gamma_k]
\]
defines a descending filtration on $\pi$. We denote the mapping class group of $\Sigma$ relative to the boundary by $\mathcal{M}$. Note that any group homomorphism $\varphi\colon \mathcal M \to \mathcal M$ preserves this filtration: $\varphi(\Gamma_k) \subset \Gamma_k$. This filtration induces the \textit{Johnson filtration} on $\mathcal M$: for $k\geq 1$, we set
\[
	\mathcal{M}(k) = \Ker(\mathcal M \to \Aut_\mathrm{grp}(\pi/\Gamma_{k+1}))\,.
\]
The first term $\mathcal{M}(1)$ is called the \textit{Torelli group} of the surface $\Sigma$. The Johnson filtration is central: $[\mathcal{M}(k), \mathcal{M}(\ell)] \subset \mathcal{M}(k+\ell)$. Therefore, we have positively graded Lie algebras (over $\mathbb{Z})$:
\[
	\mathrm{gr}(\pi) = \bigoplus_{k\geq 1} \Gamma_k/\Gamma_{k+1}\;\mbox{ and }\;\mathrm{gr}(\mathcal{I}) = \bigoplus_{k\geq 1} \mathcal{M}(k)/ \mathcal{M}(k+1)\,.
\]
The former is canonically isomorphic as a graded Lie algebra to the free Lie algebra $L(H_\mathbb{Z})$ over $H_\mathbb{Z}$. The element of $\pi$ representing $\partial\Sigma$ lives in $\Gamma_2$ and we denote its class in $\mathrm{gr}(\pi)$ by $\omega\in \Gamma_2/\Gamma_3$.

Denoting by $\Der_\omega^+(\mathrm{gr}(\pi))$ the Lie algebra of all derivations of $\mathrm{gr}(\pi)$ with positive degree and annihilating the element $\omega$, the \textit{Johnson homomorphism} $\tau_\mathbb{Z}$ (in the formulation of Morita \cite{morita}) is the Lie algebra map defined by 
\begin{align*}
	\tau_\mathbb{Z} \colon \mathrm{gr}(\mathcal{I})&\to \Der_\omega^+(\mathrm{gr}(\pi)) \cong \Der_\omega^+(L(H_\mathbb{Z}))\\
	[\varphi] &\mapsto ([\gamma] \mapsto [\varphi(\gamma)\gamma^{-1}])
\end{align*}
for $\gamma\in\pi$ and $\varphi\in\mathcal{M}$. We put $H = H_1(\Sigma;\mathbb{K})$ and denote the $\mathbb{K}$-linear extension of $\tau_\mathbb{Z}$ by
\[
	\tau\colon \mathrm{gr}(\mathcal{I})\otimes_\mathbb{Z}\mathbb{K} \to \Der_\omega^+(L(H))
\]
where $L(H)$ is the free $\mathbb{K}$-Lie algebra over $H$. We can check that $\tau_\mathbb{Z}$ and $\tau$ are injective, so the main problem is to study their image (Morita's problem, see Problem 32 of \cite{morita2}, for example). It is known to be not surjective, and the image of $\tau$ is described by the result of Hain (\cite{hain1}, \cite{hain2}): for $g\geq 3$, it is generated by the degree-$1$ part of $\Der_\omega^+(L(H))$. Therefore, the remaining problem (at least for $g\geq 3$) is to determine the defining system of equations of the image.

In \cite{morita}, Morita introduced a non-trivial map $\Tr^\mathsf{M}\colon \Der_\omega^+(L(H)) \to S(H)$ where $S(H)$ is the symmetric algebra over $H$, such that $\Tr^\mathsf{M}\circ\,\tau = 0$ on degree $\geq 2$. Later, Enomoto and Satoh refined Morita's map to the \textit{Enomoto--Satoh (ES) trace}
\[\
	\Tr^\mathsf{ES}\colon \Der^+_\omega(L(H)) \hookrightarrow \Der(T(H))^{(\geq 0)} \cong H^*\otimes T(H)^{(\geq 1)} \xrightarrow{\mathrm{cont}} T(H) \twoheadrightarrow |T(H)|\,,
\]
where $T(H)$ is the tensor algebra over $H$, $|T(H)| = T(H)/[T(H), T(H)]$ is the \textit{trace space} (i.e., the cyclic quotient) of $T(H)$, the superscript ${}^{(\geq d)}$ indicates the degree $\geq d$ part and the contraction is taken with the first component in $T(H)^{(\geq 1)}$. The important property of the ES trace is that it is a Lie algebra $1$-cocycle: we have $\Tr^\mathsf{ES}\circ\,\tau = 0$ on degree $\geq 2$ by the result of Enomoto and Satoh \cite{es}, whose proof is based on Satoh's work \cite{satoh}, which essentially relies on the $1$-cocycle property. We know that $\Imag\tau\subset \Ker\Tr^\mathsf{ES}$ is still a proper inclusion in degree $\geq 6$; in fact, the cokernel of $\tau$ is determined up to degree $9$ by combining the work of Morita--Sakasai--Suzuki as mentioned in the introduction (\cite{mss2} and their unpublished result), Garoufalidis--Getzler \cite{garget} and Kupers--Randal-Williams \cite{krw}. These are pointed out by Takuya Sakasai in personal communication. 

The ES trace is obtained from and equivalent to the graded version of the Turaev cobracket
\[
	\delta\gr\colon |T(H)|\to |T(H)|\otimes |T(H)|
\]
in the following sense: we have a canonical isomorphism of Lie algebras
\[
	|T(H)|^{(\geq 3)} \xrightarrow{\cong} \Der^+_\omega (T(H))
\]
where the domain is equipped with the graded version of the Goldman bracket $[\cdot,\cdot]\gr$ (or Schedler's Necklace Lie bracket, see \cite{necklace}) and via this isomorphism, the ES trace is obtained as the reduction $(\id\otimes\, \varepsilon)\circ\delta\gr$ of the graded Turaev cobracket by the augmentation map $\varepsilon\colon T(H)\to\mathbb{K}$. For more details, see Section 9 of \cite{akkn}.\\

\section{Ribbon Graph Complexes}\label{sec:rg} 

The above maps $[\cdot,\cdot]\gr$ and $\delta\gr$ naturally appear in the framework of ribbon graph complexes by Merkulov and Willwacher. In this section, we briefly recall some constructions and the results in their paper \cite{ribbon}.

First of all, a \textit{properad} is a generalisation of an operad, which deals with multiple outputs. (For the general theory, see \cite{properad}.) For $d\in\mathbb{Z}$, we denote by $\mathcal Lieb_d$ the properad governing the (degree-shifted) Lie bialgebra. In \cite{ribbon}, they constructed the properad $\mathcal{RG}ra_d$ of (connected, oriented and labelled) ribbon graphs together with the canonical maps of properads
\[
	\mathcal{L}ieb_d \xrightarrow{s^*} \mathcal{RG}ra_d \xrightarrow{\rho_W} \mathcal End_{|T(W)|}
\]
for any graded vector space $W$ with a degree $(1-d)$ skew-symmetric pairing $W\otimes W\to \mathbb{W}[1-d]$. In particular, the map $s^*$ is specified by
\begin{align*}
	(\textrm{Lie bracket}) \mapsto \RGbar\;\textrm{ and }\; (\textrm{Lie cobracket}) \mapsto G_1 := \RGtadpole\,.
\end{align*}
The map $\rho_W$ ensures that each ribbon graph gives a map of the form $|T(W)|^{\otimes n} \to |T(W)|^{\otimes m}$ and the map $s^*$ says that the vector space $|T(W)|$ is automatically a Lie bialgebra with degree shift. In the case of $d=1$ and $W = H$ with the intersection pairing, the structure map of this Lie bialgebra is the pair $([\cdot,\cdot]\gr,\,\delta\gr)$ by the construction of the map $\rho_H$.

We have another map
\[
	\mathcal{L}ieb_d \xrightarrow{s_1^*} \mathcal{RG}ra_d
\]
of properads, which sends the generator corresponding to the Lie cobracket in $\mathcal{L}ieb_d$ to zero instead of $G_1$. The machinery called the \textit{deformation complex} (see \cite{merval} for a nice exposition) yields the \textit{ribbon graph complex}
\begin{align*}
	(\mathsf{RGC}_d, \delta) := \mathsf{Def}(\mathcal{L}ieb_d \xrightarrow{s_1^*} \mathcal{RG}ra_d)\,.
\end{align*}
For a ribbon graph $G$, we denote by $V(G)$ the set of vertices, $E(G)$ the set of edges and $B(G)$ the set of boundary components. As a vector space, $\mathsf{RGC}_d$ is spanned by ribbon graphs with labellings on $V(G)$, $E(G)$ and $B(G)$ and the choice of direction of edges, and we can permute labellings and flip edges with signs depending on the integer $d$. The differential $\delta$ can be read off from the morphism $s_1^*$, and in this case, the differential is given by the vertex expansion. For a ribbon graph $G$, its degree $|G|_d$ in $\mathsf{RGC}_d$ is given by 
\[
	|G|_d = d(\#V(G) + \#B(G) - 2) + (1-d) \#E(G)\,,
\]
which is naturally obtained from the definition of the deformation complex (see Section 4.3 of \cite{ribbon}). The complexes $\mathsf{RGC}_d$ are isomorphic to each other up to degree shift, and the isomorphism is given by the identity map of the underlying vector space.

By the result of Strebel, Harer, Penner, Kontsevich and many others (see \cite{mp}, \cite{ribbon} and references there), the cohomology in the case of $d=0$ is described as
\begin{align}
	H^k(\mathsf{RGC}_0,\delta) \cong \prod_{\substack{g\geq 0, n\geq 1,\\2-2g-n<0}} H_{6g-6+3n-k}(\mathbf{M}_g^n;\mathbb{K})\otimes_{S_n}\mathrm{sgn}_n \oplus \left\{\begin{aligned}
		&\mathbb{K}& & \mbox{if } k\equiv 1 \,\mathrm{mod}\,4,\\
		& 0 & & \mbox{otherwise.}
	\end{aligned} \right.\label{eq:moduliisom}
\end{align}
Here, $\mathbf{M}_g^n$ denotes the moduli space of Riemann surfaces of genus $g$ with $n$ marked points, and $\mathrm{sgn}_n$ is the sign representation of the symmetric group $S_n$. This isomorphism is given in the way that the factor for the pair $(g,n)$ corresponds to ribbon graphs with $(g,n) = (g(G), \#B(G))$ where $g(G)$ is the genus of (the fattening of) $G$.

On the other hand, in the case of $d=1$, we have the canonical identification
\begin{align*}
	\mathsf{Def}(\mathcal{L}ieb_1 \xrightarrow{\rho_W\circ s_1^*} \mathcal End_{|T(W)|}) \cong C^\bullet_\mathrm{CE}(|T(W)|, \wedge|T(W)|)
\end{align*}
and hence the cochain map
\[
	\mathsf{RGC}_1 \xrightarrow{\rho_W} C^\bullet_\mathrm{CE}(|T(W)|, \wedge|T(W)|)\,.
\]
Here $\mathsf{RGC}_1$ is equipped with another grading by the number of vertices of a ribbon graph, which we call the \textit{vertex grading}. The unique graph $G_1$ with only one vertex and one edge is a $1$-cocycle in $\mathsf{RGC}_1$, and is sent to $\delta\gr$. Since the construction of the ribbon graph complex does not depend on the surface $\Sigma$, the graph $G_1$ is the ``universal'' cohomology class giving the ES trace in this perspective.

\begin{remark}
In the case of $W = H = H_1(\Sigma_{g,1};\mathbb{K})$, we can check that the cochain map $\rho_H$ is asymptotically injective as $g \to + \infty$.
\end{remark}

For a similar construction of the map of complexes for the graph homology, see the paper \cite{hamlaz} by Hamilton and Lazarev. \\

\section{The Main Result}\label{sec:mainthm}

Based on the above observations, we want to find $1$-cocycles with respect to the vertex grading in $\mathsf{RGC}_1$: if there is such a $1$-cocycle $G$, it follows that $G\circ\,\tau = 0$ by combining the construction of $s^*$ with Hain's result, giving a possible candidate for the detection of the Johnson cokernel. Unfortunately, the main result of this paper is the following:
\begin{theorem}\label{thm:nococyc}
The space of $1$-cocycles with respect to the vertex grading is $1$-dimensional and it is spanned by $G_1$.
\end{theorem}

Before the proof, recall the result of Church--Farb--Putman \cite{cfp} and Morita--Sakasai--Suzuki \cite{mss} which says that the virtual top-dimensional rational cohomology $H^{4g-3}(\mathbf{M}_g^1; \mathbb{Q})$ vanishes for $g\geq 2$.\\

\noindent\textbf{Proof of Theorem \ref{thm:nococyc}.} For a ribbon graph with $\#V(G) = 1$, we have
\[
	|G|_0 = \#E(G) = 2g(G) + \#B(G) - 1\,.
\]
and hence
\begin{align*}
	6g - 6 + 3n - |G|_0 &= (6g - 6 + 3n) - (2g + n - 1)\\
	&= 4g + 2n - 5\,.
\end{align*}
Therefore, we have
\[
	H^\bullet(\mathsf{RGC}_0,\delta)\cap (\#V(G) = 1) = \mathbb{K}G_1 \oplus \prod_{\substack{g\geq 0, n\geq 1,\\2-2g-n<0}}H_{4g+2n-5}(\mathbf{M}_g^n;\mathbb{K})\otimes_{S_n}\mathrm{sgn}_n
\]
by the isomorphism \eqref{eq:moduliisom}. 

By the result of Harer \cite{harer}, the virtual cohomological dimension of $\mathbf{M}_g^n$ with $2-2g-n<0$ is equal to 
\begin{align*}
	\mathrm{vcd}(\mathbf{M}_g^n) = \left\{\begin{aligned}
		&n - 3& &(g = 0),\\
		&4g + n - 4& &(g\geq 1).
	\end{aligned}\right.
\end{align*}
For $n\geq 2$ (and hence $n\geq 3$ when $g=0$), the factor
\[
	H_{4g+2n-5}(\mathbf{M}_g^n;\mathbb{K})\otimes_{S_n}\mathrm{sgn}_n
\]
vanishes since $4g+2n-5 > \mathrm{vcd}(\mathbf{M}_g^n)$ and we are working on a field of characteristic zero. For $n=1$, the factor is equal to $H_{4g-3}(\mathbf{M}_g^1;\mathbb{K})$, and it vanishes for $g\geq 2$ by the result mentioned above. For $(g,n) = (1,1)$, the moduli space $\mathbf{M}_1^1$ is the modular curve, which is a cusped disk, so we have $H_1(\mathbf{M}_1^1;\mathbb{K}) = 0$.

Finally, we have the isomorphism
\[
	H^\bullet(\mathsf{RGC}_0,\delta)\cap (\#V(G) = 1) \cong H^\bullet(\mathsf{RGC}_1,\delta)\cap (\#V(G) = 1)
\]
since the complexes $\mathsf{RGC}_d$ are isomorphic via the identity map. In addition, cohomology classes with $\#V(G) = 1$ are equivalent to $1$-cocycles since the differential increases $\#V(G)$ by $1$ and there are no ribbon graphs with $\#V(G) = 0$. This completes the proof.\qed\\

Still, we can ask if there is a higher cocycle in $\mathsf{RGC}_1$ which detects a part of the Johnson cokernel that the ES trace (or even the Conant--Kasabov--Vogtmann trace in \cite{ckv2}) does not. Since the cohomology of $\mathsf{RGC}_1$ is infinite-dimensional, we have many cohomology classes (let alone higher cocycles; they are abundant). In general, if we have a $k$-cocycle $f$ over a Lie algebra $\mathfrak g$, we can define the Lie subalgebra of annihilators
\[
	\mathrm{Ann}(f) = \{x\in\mathfrak g\colon f(x, y_2,\dotsc y_k) = 0 \mbox{ for any } y_2,\dotsc y_k\in \mathfrak g\},
\]
which is a generalisation of the kernel of a $1$-cocycle. The precise statement of the above question would be:
\begin{question}
For $k>1$, is there any $k$-cocycle $c$ of $\mathsf{RGC}_1$ with respect to the vertex grading such that $\Imag \tau \subset \mathrm{Ann}(\rho_H(c))$ but $\mathrm{Ann}(G_1) \not\subset \mathrm{Ann}(\rho_H(c))$ for large $g$?
\end{question}

\begin{example}
For $k\equiv 1\,\mathrm{mod}\,4$, we have the $k$-cocycle given by the unique bivalent ribbon graph $G_k$ with $k$ vertices and $k$ edges. This spans the one-dimensional subspace of $H^k(\mathsf{RGC}_0,\delta)$, which is the second summand of \eqref{eq:moduliisom}. However, we can check that $\Imag \tau \not\subset \mathrm{Ann}(\rho_H(c))$.
\end{example}
The graph $G_k$ is used in the construction of an algebraic extension of the Turaev cobracket by the author in \cite{toyo3} motivated by the relation between the Turaev cobracket and the ES trace.\\

We have considered objects in the associative world so far, but since $\Der(L(H))$ lives in the Lie world, we ask the following:
\begin{question}\label{q:Lie}
Can we do a similar discussion to the above using the Lie graph complex? For example, can we endow the space of Lie graphs with a structure similar to a properad?
\end{question}

As mentioned in the introduction, the construction of the trace map by Conant, Kassabov and Vogtmann naturally utilises Lie graphs. They showed that it takes its value in the cohomology $H^{2n-3}(\mathrm{Out}(F_n); \overline{T(V)^{\otimes n}})$ (see \cite{ck}) based on their work on hairy graph complexes, which include the (hairy) Lie graph complex. It would be interesting if we could recover their construction using the Lie graph complex and further extend their trace map.\\

\small
\bibliographystyle{alphaurl}
\bibliography{jrgc.bib}

\end{document}